# Feedback Control of Scalar Conservation Laws with Application to Density Control in Freeways by Means of Variable Speed Limits

## Iasson Karafyllis[*] and Markos Papageorgiou[**]

[*]Dept. of Mathematics, National Technical University of Athens,
Zografou Campus, 15780, Athens, Greece, email: iasonkar@central.ntua.gr

[**] Dynamic Systems and Simulation Laboratory, Technical University of Crete,
Chania, 73100, Greece, email: markos@dssl.tuc.gr

**Abstract**

The paper provides results for the stabilization of a spatially uniform equilibrium profile for a scalar conservation law that arises in the study of traffic dynamics under variable speed limit control. Two different control problems are studied: the problem with free speed limits at the inlet and the problem with no speed limits at the inlet. Explicit formulas are provided for respective feedback laws that guarantee stabilization of the desired equilibrium profile. For the first problem, global asymptotic stabilization is achieved; while for the second problem, regional exponential stabilization is achieved. Moreover, the solutions for the corresponding closed-loop systems are guaranteed to be classical solutions, i.e., there are no shocks. The obtained results are illustrated by means of a numerical example.

**Keywords:** Conservation laws, hyperbolic PDE systems, feedback stabilization.

## 1. Introduction

Scalar conservation laws have been studied extensively in the last two decades; see for instance [1, 10, 19]. Control problems related to 1-D scalar conservation laws have been investigated in many works; see [1, 2, 7, 8, 11, 17, 20, 24, 25]. Traffic control problems related to one or two 1-D conservation laws have recently been studied in [18, 27, 28, 29, 30].

In this work, we study a specific 1-D scalar conservation law, which is described by the following first-order Partial Differential Equation (PDE):

$$\frac{\partial \rho}{\partial t}(t,x) + \frac{\partial}{\partial x}\big(u(t,x)f(\rho(t,x))\big) = 0 \tag{1}$$

where $\rho$ is the state variable, $u$ is the control input, $t \geq 0$ denotes time and $x \geq 0$ denotes spatial position. To the best of our knowledge, this is the first work that studies control problems related to the PDE (1). This particular 1-D scalar conservation law arises in the study of traffic dynamics in a freeway, where Variable Speed Limits (VSL) can be applied continuously (both in time and space) along the freeway. The PDE (1) is a variation of the standard first-order LWR model (see [21, 26]) and was used in [12] for traffic networks with VSL which depend only on time. In the context of traffic dynamics, the state is the vehicle density $\rho$ (in veh/km), and the control input $u$ is the speed limit ratio ($u=1$ is the case where no speed limits are applied, and $u=0$ is the case where no vehicle movement is allowed). Thus, the PDE (1) is accompanied by the constraint $0 < u \leq 1$.

It should be noted at this point that, although the PDE (1) has not been studied so far for control purposes, the effect of VSL on traffic flow and the exploitation of this effect for improved traffic flow efficiency have been studied mainly by use of discretized models (see [1, 3, 4, 5, 6, 14, 15, 22, 23]).

The objective of the present work is the stabilization of a spatially uniform equilibrium profile by means of a feedback law, and we consider two different variants of the related control problem:
1) The control problem where we are free to use any speed limit at the inlet ($x = 0$); this implies that the inlet itself is considered controllable.
2) The control problem where we are not allowed to have speed limits at the inlet, i.e., $u(t,0) \equiv 1$; thus, the inlet for this problem is not controllable, but internal control is still allowed.

The motivation for the study of the second control problem stems from the concern to avoid the creation of queues at the entrance of the freeway, which may be required in specific applications. The two control problems are qualitatively different, with the second problem being more demanding than the first one. The qualitative difference is mathematically expressed by the imposition of a boundary condition in the second control problem.

We provide explicit formulas for the feedback laws that guarantee exponential stabilization of the desired equilibrium profile. Moreover, the solutions of the closed-loop systems (whose existence is established in both cases by Banach's fixed point theorem) are classical solutions, i.e., continuously differentiable; therefore, there are no shocks. However, for the second control problem, we cannot achieve global stabilization (i.e., for all physically relevant initial conditions). The proposed feedback law for the second control problem achieves regional stabilization, and the region where stabilization can be achieved depends on the parameters of the controller. This is conform with physical intuition and traffic engineering experience, which has shown that extreme congestion phenomena cannot be handled without controlling the inlet flow (via VSL-induced mainstream metering as proposed in [4, 6] or via ramp metering or both).

The structure of the present work is as follows: Section 2 is devoted to the presentation of the model and the two control problems related to the 1-D scalar conservation law (1). Section 3 provides the statements of the main results (Theorem 3.1 and Theorem 3.2), as well as a discussion thereof. The proofs of the main results are provided in Section 4. A simple illustrative example is presented in Section 5. Concluding remarks are provided in Section 6.

**Notation.** Throughout this paper, we adopt the following notation.

* $\Re_+ := [0, +\infty)$. Let $U \subseteq \Re^n$ be a set with non-empty interior and let $\Omega \subseteq \Re$ be a set. By $C^0(U; \Omega)$, we denote the class of continuous mappings on $U$, which take values in $\Omega$. By $C^k(U; \Omega)$, where $k \geq 1$, we denote the class of continuous functions on $U$, which have continuous derivatives of order $k$ on $U$ and take values in $\Omega$.

* Let $\rho : I \times [0, L] \to \Re$ be given, where $I \subseteq \Re$ is an interval and $L > 0$ is a constant. We use the notation $\rho[t]$ to denote the profile of $\rho$ at certain $t \in I$, i.e., $(\rho[t])(z) = \rho(t, z)$ for all $z \in [0, L]$.

* Let $\rho \in C^0([0, L]; \Re)$ be given, where $L > 0$ is a constant. For every $\rho^* \in \Re$ we use the notation $\|\rho - \rho^*\|_\infty$ to denote the sup norm of the function $f(x) = \rho(x) - \rho^*$, i.e., $\|\rho - \rho^*\|_\infty := \max_{0 \leq x \leq L} \left(|\rho(x) - \rho^*|\right)$.

## 2. Description of two control problems

*2.1 The LWR model with VSL*

VSL affect the fundamental diagram (flow-density curve) in a nonlinear way (see Chapter 2 in [6] and also [9, 13, 16]). In general, for every freeway there exists a nonlinear function $F \in C^2([0, \rho_{max}] \times [0,1]; [0, q_{max}])$, where $\rho_{max} > 0$ is the physical upper bound of density in the particular road and $q_{max} > 0$ is the upper bound (capacity) of the flow of vehicles in the road, with $F(0,l) = F(\rho,0) = 0$ for all $(\rho,l) \in [0, \rho_{max}] \times [0,1]$ and $F(\rho,l) > 0$ for all $(\rho,l) \in (0, \rho_{max}] \times (0,1]$, so that the flow $q$ of vehicles at a point in freeway with density $\rho \in (0, \rho_{max}]$ is given by the equation

$$q = F(\rho, l) \tag{2}$$

where $l \in (0,1]$ is the ratio of the imposed speed limit divided by the maximum mean speed that can be exhibited in the freeway without speed limit.

Assuming that the speed limit ratio, i.e., $l$, can be manipulated continuously (both in time and space), i.e., assuming that $l = l(t,x)$, we obtain the modified LWR model with VSL for $(t,x) \in \Re_+ \times [0,L]$:

$$\frac{\partial \rho}{\partial t}(t,x) + \frac{\partial}{\partial x}\big(F(\rho(t,x), l(t,x))\big) = 0 \tag{3}$$

where $L > 0$ is the length of the freeway, $\rho(t,x)$ is the vehicle density, $t \geq 0$ is time and $x \in [0,L]$ is the spatial variable (position).

We next use the following assumption.

**(H1)** *There exists a continuous function $\tilde{l} : (0, \rho_{max}] \to (0,1]$ and a constant $\delta \in (0, \rho_{max}]$ such that*

$$\frac{\partial F}{\partial l}(\rho, l) > 0 \text{, for all } l \in (0, \tilde{l}(\rho)) \text{ and } \rho \in (0, \rho_{max}] \tag{4}$$

$$F(\rho, \tilde{l}(\rho)) = F(\rho, 1) \text{, for all } \rho \in (0, \rho_{max}] \tag{5}$$

$$\tilde{l}(\rho) = 1, \text{ for all } \rho \in (0, \delta] \tag{6}$$

This assumption reflects the physical effect of VSL on the flow-density curve in a rather general way: the application of any speed limit $l < 1$ reduces the flow for $\rho \in (0, \delta]$ and has an arbitrary effect on the flow for $\rho \in (\delta, \rho_{max}]$.

**Remark 2.1:** It should be noted that Assumption (H1) holds for the class of functions

$$F(\rho, l) = A\rho l \exp\left(-\frac{1}{\gamma}\left(\frac{b\rho}{1+a-al}\right)^\gamma\right) \tag{7}$$

with $A, b, \gamma > 0$, $a \geq 0$, which is a special case of the class of functions proposed in [6] for the description of the effect of VSL to the flow-density curve. More specifically, using $\frac{1}{A\rho}\frac{\partial F}{\partial l}(\rho,l)\exp\left(\frac{1}{\gamma}\left(\frac{b\rho}{1+a-al}\right)^{\gamma}\right)=1-\frac{al\,(b\rho)^{\gamma}}{(1+a-al)^{1+\gamma}}$, we can guarantee that $\delta=\rho_{\max}$ when $a\,(b\rho_{\max})^{\gamma}\leq 1$ and $\delta=b^{-1}a^{-1/\gamma}<\rho_{\max}$ when $a\,(b\rho_{\max})^{\gamma}>1$. The function $\tilde{l}:(0,\rho_{\max}]\to(0,1]$ for the case (7) is defined for each $\rho\in(0,\rho_{\max}]$ as the smallest solution $l$ of the equation $(1+a-al)^{\gamma}\left(1+\gamma(b\rho)^{-\gamma}\ln(l)\right)=1$.

Assumption (H1) allows us to define a continuous function $g:\bigcup_{\rho\in(0,\rho_{\max})}(\{\rho\}\times(0,F(\rho,1)])\to(0,1]$ so that the following equations hold for all $\rho\in(0,\rho_{\max}]$ and $y\in(0,F(\rho,1)]$:

$$0<g(\rho,y)\leq\tilde{l}(\rho),\ y=F(\rho,g(\rho,y)) \tag{8}$$

Define:

$$f(\rho):=F(\rho,1),\ \text{for}\ \rho\in[0,\rho_{\max}] \tag{9}$$

Under Assumption (H1) and definition (9), for each $u\in(0,1]$ and $\rho\in(0,\rho_{\max}]$, we are in a position to find $l\in(0,\tilde{l}(\rho)]$ so that $uf(\rho)=F(\rho,l)$ (namely, $l=g(\rho,uf(\rho))$). Therefore, we may consider the simplified LWR model (1) with VSL for $(t,x)\in\Re_{+}\times[0,L]$. Traffic flow theory allows us to use the following assumption for the function $f\in C^{2}([0,\rho_{\max}];[0,q_{\max}])$, defined by (9).

**(H2)** *The function* $f\in C^{2}([0,\rho_{\max}];[0,q_{\max}])$, *where* $\rho_{\max}>0$, $q_{\max}>0$ *is a function for which there exists* $\rho_{cr}\in(0,\rho_{\max})$ *with the following properties: (i)* $f'(\rho)>0$ *for all* $\rho\in[0,\rho_{cr})$, $f'(\rho_{cr})=0$, *(ii)* $f''(\rho)<0$ *for all* $\rho\in[0,\rho_{\max}]$, *and (iii)* $f(0)=0$ *and* $f(\rho)>0$ *for all* $\rho\in(0,\rho_{\max}]$.

**Remark 2.2:** Assumption (H2) introduces the critical density $\rho_{cr}\in(0,\rho_{\max})$, which in traffic engineering is the density value that produces the maximum flow (capacity).

*2.2 Two Control Problems*

Let $\rho^{*}\in(0,\delta)$ be the given set point for density.

1st Problem: Free Speed Limit at Inlet

The control objective is to construct a feedback law of the form

$$u(t,x)=K(\rho[t],x),\ \text{for}\ (t,x)\in\Re_{+}\times[0,L] \tag{10}$$

with

$$K(\rho[t],x)\in(0,1],\ \text{for}\ (t,x)\in\Re_{+}\times[0,L] \tag{11}$$

so that for any initial condition $\rho_{0}\in C^{1}([0,L];(0,\rho_{\max}])$, the solution $\rho[t]$ of the closed-loop system (1) with (10) and initial condition $\rho[0]=\rho_{0}$ exists for all $t\geq 0$, is unique and satisfies

$$\lim_{t\to+\infty}(\rho(t,x))=\rho^{*},\ \lim_{t\to+\infty}(u(t,x))=1,\ \text{for}\ x\in[0,L]. \tag{12}$$

The requirement $\lim_{t\to+\infty}(u(t,x))=1$ in (12) is important, because it implies that the VSL may be employed to address any appearing congestion problem, hence the freeway will practically operate without speed limits after an initial transient period, i.e. after the problem has been tackled; notice that this is a consequence of (4), (5), (6) and (9) (recall that $\rho^* \in (0,\delta)$, where $\delta > 0$ is the constant involved in (H1)).

2nd Problem: No Speed Limit at Inlet

One possible issue with the 1st problem is that the inlet flow, which is equal to $u(t,0)f(\rho(t,0))$, may become small for a transient period. This may cause the creation of queues at the entrance of the freeway (or at upstream on-ramps, in case of ramp metering). One way to avoid the creation of queues is to require that no speed limit is applied at the entrance of the freeway, i.e.,

$$u(t,0)=1, \text{ for } t \geq 0 \tag{13}$$

Moreover, in this case, the flow at the entrance of the freeway is assumed to be equal to the nominal flow, i.e.,

$$\rho(t,0)=\rho^*, \text{ for } t \geq 0 \tag{14}$$

The control objective is to construct a feedback law of the form (10), satisfying (11) and $K(\rho[t],0)=1$, so that for a large number of initial conditions $\rho_0 \in C^1([0,L];(0,\rho_{max}])$ with $\rho_0(0)=\rho^*$, the solution $\rho[t]$ of the closed-loop system (1), (14) with (10) and initial condition $\rho[0]=\rho_0$ exists for all $t \geq 0$, is unique and satisfies (12).

The 2nd control problem described above is not expected to be solvable for arbitrary set points $\rho^* \in (0,\rho_{cr})$ and arbitrary initial conditions $\rho_0 \in C^1([0,L];(0,\rho_{max}])$. In fact, physical intuition and traffic engineering experience indicate that extreme congestion phenomena (e.g., $\rho_0(x)=\rho_{max}$ for some $x \in (0,L]$) may be hard to handle without controlling the inlet flow.

## 3. Main Results

*3.1 Two possible solutions of the control problems*

The following theorem guarantees that the 1st problem is globally solvable, i.e., for all physically relevant initial conditions.

**Theorem 3.1 (Global Asymptotic Stabilization in the sup norm):** *Suppose that Assumption (H2) holds. Let $\rho^* \in (0,\delta)$ and let $k \in (0,1/(L\rho^*))$ be given constants. Define for all $x \in [0,L]$ and $w \in C^0([0,L];(0,\rho_{max}])$*

$$M(w,x)=1/\left(1+k\int_0^x (w(s)-\rho^*)ds\right). \tag{15}$$

*Then there exists a non-decreasing function $c:(0,\rho_{max}] \to (0,q_{max}]$ such that for every $\rho_0 \in C^1([0,L];(0,\rho_{max}])$ the initial value problem (1) with $\rho[0]=\rho_0$ and*

$$u(t,x) = \frac{\min_{z\in[0,L]} \left( f(\rho(t,z))M(\rho[t],z) \right)}{f(\rho(t,x))M(\rho[t],x)}, \text{ for } t\geq 0,\ x\in[0,L] \tag{16}$$

*has a unique solution* $\rho \in C^1(\Re_+ \times [0,L];(0,\rho_{\max}])$, *which satisfies the following estimates:*

$$\|\rho[t]-\rho^*\|_\infty \leq \exp(-ct)\|\rho_0-\rho^*\|_\infty, \text{ for all } t\geq 0 \tag{17}$$

$$\lim_{t\to+\infty} u(t,x) = 1, \text{ for all } x\in[0,L] \tag{18}$$

*where* $c = c\left( \min_{z\in[0,L]}(\rho_0(z)) \right)$.

The following result guarantees that the 2nd control problem is solvable for a certain class of initial conditions (regional stabilization).

**Theorem 3.2 (Regional Exponential Stabilization in the sup norm):** *Suppose that Assumption (H2) holds. Let* $\rho^* \in (0,\min(\delta,\rho_{cr},\rho_{\max}/2))$ *be a given constant. Then there exist constants* $\sigma,\gamma,c>0$ *such that for every* $\rho_0 \in C^1([0,L];(0,\rho_{\max}])$ *with* $\rho_0(0)=\rho^*$ *and*

$$f(\rho^*) + \sigma\int_0^x (\rho_0(s)-\rho^*)ds - \gamma\frac{x^2}{2}\|\rho_0-\rho^*\|_\infty \leq f(\rho_0(x)), \text{ for } x\in[0,L] \tag{19}$$

*the initial-boundary value problem (1) with (14),* $\rho[0]=\rho_0$ *and*

$$u(t,x) = (f(\rho(t,x)))^{-1}\left( f(\rho^*) + \sigma\int_0^x (\rho(t,s)-\rho^*)ds - \gamma\frac{x^2}{2}\|\rho[t]-\rho^*\|_\infty \right), \text{ for } t\geq 0,\ x\in[0,L] \tag{20}$$

*has a unique solution* $\rho \in C^1(\Re_+ \times [0,L];(0,\rho_{\max}])$, *which satisfies estimates (17), (18) as well as the following estimate:*

$$0 < u(t,x) \leq 1, \text{ for } t\geq 0,\ x\in[0,L] \tag{21}$$

*3.2 Discussion of Main Results*

We provide below a list of comments for the main results.

**1)** The solutions of both control problems guarantee asymptotic stability in the sup norm.

**2)** The 2nd control problem is more demanding than the 1st control problem. The reason that explains this difference between the two control problems is related to the inlet flow; while the inlet flow may be modulated in the 1st problem, it remains constant (and equal to $f(\rho^*)$, i.e. typically very high) for the 2nd problem. This limits the range of possible VSL actions; for example, the speed limits in the 2nd problem cannot become very small, e.g. in order to dissolve a possible downstream congestion, because this could result in high accumulation of vehicles farther upstream due to the high inflow.

**3)** The difficulty in the solution of the 2nd control problem is reflected in the fact that the control objective can be achieved only for a class of initial conditions, namely initial conditions for which (19) holds. Another fact that expresses the difficulty in the solution of the 2nd problem is that the set point $\rho^*$ cannot be equal to the critical density $\rho_{cr}$, while there is no such constraint for the 1st problem.

**4)** Since condition (19) depends on the controller gains $\sigma, \gamma > 0$, there is a degree of freedom which can be used for the enlargement of the allowable set of initial conditions. However, the proof of Theorem 2 shows that the gains $\sigma, \gamma > 0$ cannot be arbitrary (they must be sufficiently small) and that the gains affect the convergence rate (the proof of Theorem 3.2 shows that the constant $c$ involved in (17) is equal to $\sigma - \gamma L$).

**5)** Inequality (21) guarantees that for every $t \geq 0$ the solution $\rho[t]$ is a function of class $C^1([0,L];(0,\rho_{\max}])$ with $(\rho[t])(0) = \rho^*$ which also satisfies inequality (19) with $\rho[t]$ in place of $\rho_0$. Therefore, the state space for the 2nd control problem is the set

$$X_{\sigma,\gamma} := \left\{ \rho \in C^1([0,L];(0,\rho_{\max}]) : \rho(0) = \rho^*, f(\rho^*) + \sigma \int_0^x (\rho(s) - \rho^*) ds - \gamma \frac{x^2}{2} \|\rho - \rho^*\|_\infty \leq f(\rho(x)) \right\}.$$

In other words, for every $\rho_0 \in X_{\sigma,\gamma}$, it follows that $\rho[t] \in X_{\sigma,\gamma}$ for all $t \geq 0$.

**6)** The role of the parameter $\gamma > 0$ is crucial for the size of the set $X_{\sigma,\gamma}$: the smaller $\gamma > 0$, the smaller is the set $X_{\sigma,\gamma}$. However, as remarked earlier, we cannot allow $\gamma > 0$ to become arbitrarily large. A trade-off between the size of the set $X_{\sigma,\gamma}$ and the convergence rate is present here (the smaller $\gamma > 0$, the faster is the convergence rate).

**7)** The proof of Theorem 3.1 shows that the convergence rate depends heavily on the controller gain $k > 0$. However, notice that inequality (17) does not guarantee a uniform convergence rate for all initial conditions in the set $C^1([0,L];(0,\rho_{\max}])$.

**8)** The feedback law (16) can be interpreted physically: the position $x^*(t) \in [0,L]$ where

$$f(\rho(t,x^*(t)))M(\rho[t],x^*(t)) = \min_{z \in [0,L]} \left( f(\rho(t,z))M(\rho[t],z) \right)$$

is exactly the position that determines the position of the "bottleneck". The feedback law (16) imposes no speed limit at this specific position, i.e., $u(t,x^*(t)) = 1$ and reduces appropriately the flowrate at every other position. The position of the bottleneck depends on the weight term $M(\rho[t],x)$, defined by (15). This interpretation allows the physical explantation of condition (19): condition (19) guarantees that there exists an appropriate weight term (not of the form (15)) so that the bottleneck occurs at the inlet $x = 0$. Using the appropriate weight term, we construct the feedback law (20), which satisfies (13).

**9)** If $\rho_0(0) = \rho^*$ then the proof of Theorem 3.1 shows that the boundary condition (14) will also hold (although the boundary condition (14) is not imposed to the closed-loop system (1), (16)). However, as noted earlier, this fact does not imply that the inlet flow is constant. In this case ($\rho_0(0) = \rho^*$) the inlet flow is equal to $u(t,0)f(\rho^*)$, which is not necessarily equal to $f(\rho^*)$.

# 4. Proofs of Main Results

We first provide the proof of Theorem 3.1.

**Proof of Theorem 3.1:** Suppose that the initial value problem (1) with $\rho[0] = \rho_0$, (16) has a solution $\rho \in C^1(\Re_+ \times [0,L]; (0, \rho_{max}])$. It follows from (15), (16) that

$$f(\rho(t,0))u(t,0) = P(\rho[t]), \text{ for } t \geq 0 \tag{22}$$

where $P(w) = \min_{z \in [0,L]} (f(w(z))M(w,z))$. Combining (22) and (15) we get for $t \geq 0$, $x \in [0,L]$:

$$f(\rho(t,x))u(t,x) = P(\rho[t])\left(1 + k\int_0^x (\rho(t,s) - \rho^*)ds\right). \tag{23}$$

It follows from (23) and (1) that:

$$\frac{\partial \rho}{\partial t}(t,x) = -k(\rho(t,x) - \rho^*)P(\rho[t]), \text{ for } t \geq 0, \ x \in [0,L]. \tag{24}$$

Integrating (24) and using (22) and the initial condition $\rho[0] = \rho_0$, we get for $t \geq 0$, $x \in [0,L]$:

$$\rho(t,x) = \rho^* + (\rho_0(x) - \rho^*)\exp\left(-k\int_0^t g(s)ds\right), \tag{25}$$

where

$$g(t) = P(\rho[t]), \text{ for } t \geq 0 \tag{26}$$

We next show that for certain $T > 0$ (independent of $\rho_0$) the mapping $G_{\rho_0} : C^0([0,T]; \Re_+) \to C^0([0,T]; \Re_+)$ defined by

$$(G_{\rho_0} g)(t) = P(\rho[t]), \text{ for } t \in [0,T] \tag{27}$$

where $\rho : [0,T] \times [0,L] \to (0, \rho_{max}]$ is defined by (25), is a contraction. Indeed, notice that for every pair of functions $g_i \in C^0([0,T]; \Re_+)$, $i = 1, 2$, the functions $\rho_i : [0,T] \times [0,L] \to (0, \rho_{max}]$, $i = 1, 2$ satisfying $\rho_i(t,x) = \rho^* + (\rho_0(x) - \rho^*)\exp\left(-k\int_0^t g_i(s)ds\right)$ for $i = 1, 2$, $(t,x) \in [0,T] \times [0,L]$ also satisfy the estimates for $(t,x) \in [0,T] \times [0,L]$:

$$|\rho_1(t,x) - \rho_2(t,x)| \leq k|\rho_0(x) - \rho^*|\left|\int_0^t (g_2(s) - g_1(s))ds\right| \tag{28}$$

$$\leq kT|\rho_0(x) - \rho^*|\max_{0 \leq s \leq t}(|g_1(s) - g_2(s)|)$$

$$|P(\rho_1[t]) - P(\rho_2[t])| \leq \max_{z \in [0,L]} \left( |f(\rho_1(t,z))M(\rho_1[t],z) - f(\rho_2(t,z))M(\rho_2[t],z)| \right)$$

$$\leq \max_{z \in [0,L]} \left( M(\rho_1[t],z)|f(\rho_1(t,z)) - f(\rho_2(t,z))| \right) + \max_{z \in [0,L]} \left( f(\rho_2(t,z))|M(\rho_1[t],z) - M(\rho_2[t],z)| \right) \quad (29)$$

$$\leq \frac{q_{max}kL + L_f}{(1-kL\rho^*)^2} \max_{z \in [0,L]} \left( |\rho_1(t,z) - \rho_2(t,z)| \right)$$

where $L_f$ is the Lipschitz constant for $f$ on $[0, \rho_{max}]$. For the derivation of (29), we have used definition (15) as well as the facts that $f \in C^2([0,\rho_{max}];[0,q_{max}])$ and $\rho_i(t,x) \in [0, \rho_{max}]$ for $i = 1, 2$, $(t,x) \in [0,T] \times [0,L]$. Since $|\rho_0(x) - \rho^*| \leq \max(\rho^*, \rho_{max} - \rho^*)$, it follows from (28), (29) that for every $T > 0$ with

$$\frac{q_{max}kL + L_f}{(1-kL\rho^*)^2} kT \max(\rho^*, \rho_{max} - \rho^*) < 1$$

the mapping $G_{\rho_0} : C^0([0,T]; \Re_+) \to C^0([0,T]; \Re_+)$ defined by (25), (27) is a contraction. It follows from Banach's fixed point theorem and completeness of $C^0([0,T]; \Re_+)$ that the initial value problem (1) with $\rho[0] = \rho_0$, (16) has a unique solution $\rho \in C^1([0,T] \times [0,L]; (0, \rho_{max}])$. Since $T > 0$ is independent of $\rho_0$ the argument may be repeated on the intervals $[T, 2T], [2T, 3T], \ldots$ in order to construct a solution $\rho \in C^0(\Re_+ \times [0,L]; (0, \rho_{max}])$.

Notice that (25) implies that $\min \left( \min_{x \in [0,L]} (\rho_0(x)), \rho^* \right) \leq \rho(t,x) \leq \rho_{max}$ for all $(t,x) \in \Re_+ \times [0,L]$. Therefore, definitions (15), (26) and the previous inequality imply that

$$g(t) \geq \frac{\min \left\{ f(\rho) : \min \left( \min_{x \in [0,L]} (\rho_0(x)), \rho^* \right) \leq \rho \leq \rho_{max} \right\}}{1 + kL(\rho_{max} - \rho^*)} \quad (30)$$

Estimate (17) is a direct consequence of (30) in conjunction with (25) and the definition $c(s) := k \frac{\min \{ f(\rho) : \min(s, \rho^*) \leq \rho \leq \rho_{max} \}}{1 + kL(\rho_{max} - \rho^*)}$ for $s > 0$. Finally, (18) follows from (16), (17) and continuity of the functional $\min_{z \in [0,L]} (f(\rho)M(\rho,z))/M(\rho,x)$ for $\rho \in C^0([0,L]; (0, \rho_{max}])$. The proof is complete. ◁

We next provide the proof of Theorem 3.2.

**Proof of Theorem 3.2:** Let $\sigma > 0$ be a constant sufficiently small, so that

$$f(\rho^*) > \frac{\sigma L}{2}(\rho_{max} + \rho^*), \quad f'(\rho^*) > \sigma L \quad (31)$$

By virtue of Assumption (H2) there exists a unique $a \in (0, \rho_{cr} - \rho^*)$ so that

$$f'(\rho^* + a) = \sigma L \tag{32}$$

Let $\gamma > 0$ be a constant sufficiently small, so that

$$\frac{\sigma \min\{-f''(\rho): \rho \in [0, \rho_{max}]\} a^2}{2(q + \sigma L)(\rho_{max} - \rho^*)} > \gamma L \tag{33}$$

$$\sigma > \gamma L \tag{34}$$

where $q := \max\left(0, \max\{-f'(\rho): \rho \in [\rho^* + a, \rho_{max}]\}\right) > 0$. Every solution $\rho \in C^1(\Re_+ \times [0, L]; (0, \rho_{max}))$ of the initial-boundary value problem (1) with (14), $\rho[0] = \rho_0$ and (20) satisfies the following equation for $t \geq 0$, $x \in [0, L]$:

$$\frac{\partial \rho}{\partial t}(t, x) = -\sigma\left(\rho(t, x) - \rho^*\right) + \gamma x \|\rho[t] - \rho^*\|_\infty \tag{35}$$

It follows from (35) that the following equations hold for $t \geq 0$, $x \in [0, L]$:

$$\rho(t, x) = \rho^* + \exp(-\sigma t)\left(\rho_0(x) - \rho^*\right) + \gamma x \int_0^t \exp(-\sigma(t-s)) g(s) ds \tag{36}$$

$$g(t) = \|\rho[t] - \rho^*\|_\infty = \exp(-\sigma t) \max_{0 \leq x \leq L}\left(\left|\rho_0(x) - \rho^* + \gamma x \int_0^t \exp(\sigma s) g(s) ds\right|\right) \tag{37}$$

We next show that for each $T > 0$ and $\rho_0 \in C^1([0, L]; (0, \rho_{max}))$ with $\rho_0(0) = \rho^*$, the mapping $G_{\rho_0}: C^0([0, T]; \Re_+) \to C^0([0, T]; \Re_+)$ defined by

$$(G_{\rho_0} g)(t) = \exp(-\sigma t) \max_{0 \leq x \leq L}\left(\left|\rho_0(x) - \rho^* + \gamma x \int_0^t \exp(\sigma s) g(s) ds\right|\right), \text{ for } t \in [0, T] \tag{38}$$

is a contraction. Elementary manipulations allow us to show that for every pair of functions $g_i \in C^0([0, T]; \Re_+)$, $i = 1, 2$, the following inequality holds for all $t \geq 0$:

$$\left|(G_{\rho_0} g_1)(t) - (G_{\rho_0} g_2)(t)\right| \leq \frac{\gamma L}{\sigma} \max_{0 \leq s \leq t}\left(|g_1(s) - g_2(s)|\right) \tag{39}$$

Thus, since (34) holds, it follows that the mapping $G_{\rho_0}: C^0([0, T]; \Re_+) \to C^0([0, T]; \Re_+)$ defined by (38) is a contraction. It follows from Banach's fixed point theorem and completeness of $C^0([0, T]; \Re_+)$ that for each $T > 0$, the initial value problem (1) with (14), $\rho[0] = \rho_0$, (20) has a unique solution $\rho \in C^1([0, T] \times [0, L]; (0, +\infty))$. Moreover, we get from (37) for $y(t) := g(t) \exp(\sigma t)$:

$$y(t) \leq \|\rho_0 - \rho^*\|_\infty + \gamma L \int_0^t y(s) ds \tag{40}$$

Using (40) and applying Gronwall-Bellman Lemma for the continuous function $y(t) := g(t)\exp(\sigma t)$, we obtain (17) with $c := \sigma - \gamma L > 0$ (recall (34)).

Next we show that $\rho(t,x) \leq \rho_{max}$ for $t \geq 0$, $x \in [0,L]$. Suppose on the contrary that there exist $t \geq 0$, $x \in [0,L]$ such that $\rho(t,x) > \rho_{max}$. It follows that $\|\rho[t] - \rho^*\|_\infty > \rho_{max} - \rho^*$ and using (17) we get $\|\rho_0 - \rho^*\|_\infty > \rho_{max} - \rho^*$. Since $\rho_0 \in C^1([0,L];(0,\rho_{max}])$ and since $\rho_{max} > \rho^*/2$ we obtain that $\|\rho_0 - \rho^*\|_\infty \leq \rho_{max} - \rho^*$, a contradiction. Thus $\rho \in C^1(\Re_+ \times [0,L];(0,\rho_{max}])$.

The rest of the proof is devoted to showing (21). Indeed, it follows from (20), (31), (34) and since $0 < \rho(t,x) \leq \rho_{max}$ for $t \geq 0$, $x \in [0,L]$ that the inequality $u(t,x) > 0$ holds for $t \geq 0$, $x \in [0,L]$. Using (37) and the inequality

$$g(t+h) \geq \exp(-\sigma h)g(t) - \gamma L \exp(-\sigma h)\exp(-\sigma t)\int_t^{t+h}\exp(\sigma s)g(s)ds$$

we conclude that

$$\liminf_{h \to 0^+} h^{-1}(g(t+h) - g(t)) \geq -(\sigma + \gamma L)g(t), \text{ for all } t \geq 0 \quad (41)$$

We next show that $u(t,x) \leq 1$ holds for $t \geq 0$, $x \in [0,L]$, or equivalently that $f(\rho(t,x)) \geq f(\rho^*) + \sigma \int_0^x (\rho(t,s) - \rho^*)ds - \gamma \frac{x^2}{2}\|\rho[t] - \rho^*\|_\infty$ for $t \geq 0$, $x \in [0,L]$. Suppose the contrary, i.e., that there exists $t > 0$, $x \in (0,L]$ such that $f(\rho(t,x)) < f(\rho^*) + \sigma \int_0^x (\rho(t,s) - \rho^*)ds - \gamma \frac{x^2}{2}\|\rho[t] - \rho^*\|_\infty$. Define for $l \in [0,t]$:

$$\mu(l) := f(\rho(l,x)) - f(\rho^*) - \sigma \int_0^x (\rho(l,s) - \rho^*)ds + \gamma \frac{x^2}{2}\|\rho[l] - \rho^*\|_\infty \quad (42)$$

Notice that (19) implies that $\mu(0) \geq 0$ and by assumption it holds that $\mu(t) < 0$. It follows (by continuity of $\mu$) that there exists $T \in [0,t)$ with $\mu(T) = 0$ and $\mu(l) < 0$ for $l \in (T,t]$ ($T := \sup\{l \in [0,t]: \mu(l) = 0\}$). Consequently, we must have $\limsup_{h \to 0^+} h^{-1}(\mu(T+h) - \mu(T)) \leq 0$. Using (35), (37) and definition (42) we obtain for $l \in [0,t]$:

$$\limsup_{h \to 0^+} h^{-1}(\mu(l+h) - \mu(l)) \geq \sigma\left(f(\rho(l,x)) - f(\rho^*) - f'(\rho(l,x))(\rho(l,x) - \rho^*)\right)$$

$$-\sigma\mu(l) + \gamma x f'(\rho(l,x))\|\rho[l] - \rho^*\|_\infty + \gamma \frac{x^2}{2}\liminf_{h \to 0^+} h^{-1}(g(l+h) - g(l))$$

Using (37), (41), the fact that $\mu(T) = 0$ and the above inequality, we get:

$$\limsup_{h \to 0^+} h^{-1}(\mu(T+h) - \mu(T)) \geq \sigma\left(f(\rho(T,x)) - f(\rho^*) - f'(\rho(T,x))(\rho(T,x) - \rho^*)\right)$$

$$+ \gamma x \left( f'(\rho(T,x)) - (\sigma + \gamma L)\frac{x}{2} \right) \|\rho[T] - \rho^*\|_\infty$$

Finally, using the fact that $\limsup_{h \to 0^+} h^{-1}(\mu(T+h) - \mu(T)) \leq 0$ and (34) we get:

$$\sigma\left(f(\rho(T,x)) - f(\rho^*) - f'(\rho(T,x))(\rho(T,x) - \rho^*)\right) + \gamma x \left(f'(\rho(T,x)) - \sigma L\right)\|\rho[T] - \rho^*\|_\infty \leq 0 \quad (43)$$

We show next that (43) cannot hold. Notice that since $\mu(l) < 0$ for $l \in (T,t]$, it follows that $\|\rho[l] - \rho^*\|_\infty > 0$ for $l \in (T,t]$. Consequently, (17) implies that $\|\rho[T] - \rho^*\|_\infty > 0$. Moreover, notice that

$$f(\rho(T,x)) - f(\rho^*) - f'(\rho(T,x))(\rho(T,x) - \rho^*) \geq \frac{1}{2}Q(\rho(T,x) - \rho^*)^2 \quad (44)$$

where $Q := \min\{-f''(\rho) : \rho \in [0, \rho_{max}]\} > 0$ (recall Assumption (H2)). It follows from (32), (44) that (43) cannot hold for $\rho(T,x) < \rho^* + a$. For $\rho(T,x) \geq \rho^* + a$, we obtain from (43), (44) and the fact that $\|\rho[T] - \rho^*\|_\infty \leq \rho_{max} - \rho^*$:

$$\sigma Q a^2 \leq 2\gamma L(q + \sigma L)(\rho_{max} - \rho^*) \quad (45)$$

where $q := \max\left(0, \max\{-f'(\rho) : \rho \in [\rho^* + a, \rho_{max}]\}\right) > 0$.

Inequality (45) contradicts inequality (33). We conclude that (43) cannot hold. The proof is complete. ◁

## 5. Illustrative Example

Consider the case $f(\rho) = \rho \exp(-\rho)$ with

$$\rho_{max} = 1.6, \quad L = 1, \quad \rho^* = 0.7$$

Notice that in this case we have $\rho_{cr} = 1$. We apply the feedback law (16) with $k = 0.3$ and the feedback law (20) with $\sigma = 0.12$, $\gamma = 0.1$. Figures 1 and 2 show the density profiles for the initial condition $\rho_0(x) = \rho^* + 4x^2(1.2 - x)^2$, $x \in [0, L]$, while Figures 3 and 4 show the corresponding speed limit ratios $u(t,x)$.

Figures 1 and 2 show that the convergence rate of the closed-loop system (1) with (16) is faster than that of the closed-loop system (1), (14) with (20). Notice that we cannot increase the values of the controller gains $\sigma, \gamma > 0$ to speed up the convergence rate of the closed-loop system (1), (14) with (20), because, as the proof of Theorem 3.2 showed, the controller gains $\sigma, \gamma > 0$ must be sufficiently small. The speed limits are initially generally higher for the closed-loop system (1), (14) with (20) than those of the closed-loop system (1) with (16). However, the convergence rate of the speed limit ratio to 1 is faster for the closed-loop system (1) with (16) than that of the closed-loop system (1), (14) with (20).

From a physical point of view, Figures 1-4 illustrate the control mechanism leading to the dissolution of the initial congestion in the freeway. Specifically, under both control laws, the speed limit control is shown to take low values upstream of the congestion, so as to limit the respective inflows; in the first case, this involves also some inlet reduction, which is not allowed for the second case. In terms of outflow, it may be seen that the convergence rate of the density at the outlet for the closed-loop system (1), (14) with (20) is significantly slower than the convergence rate of the density at the outlet for the closed-loop system (1) with (16); this happens because the feedback law (20) cannot affect the inlet flow rate, hence it tries to maintain large outflows by keeping the density at the outlet close to the critical desnity (hence maintaining a larger outflow) for a longer period of time.

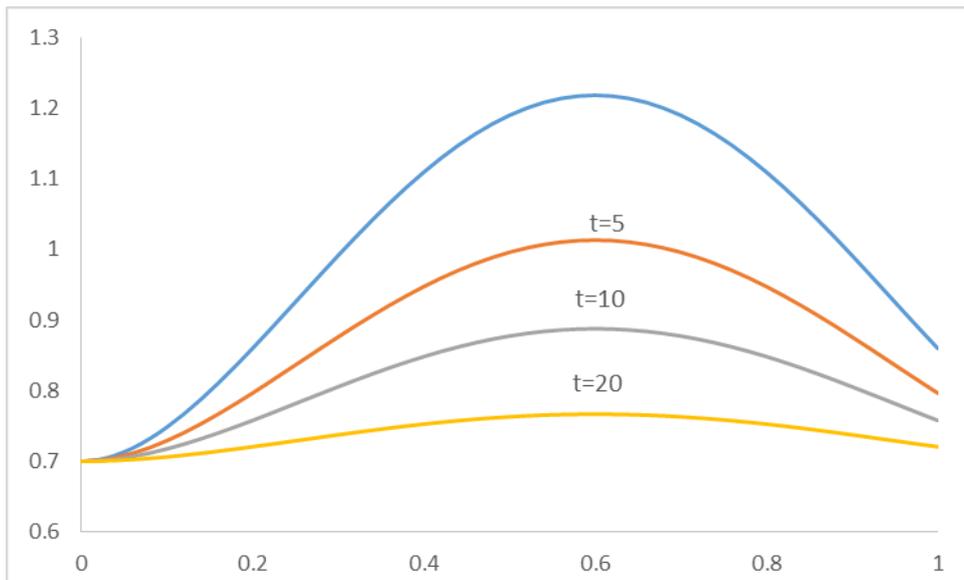

**Fig. 1:** Density ($\rho(t,x)$) profiles for the closed-loop system (1) with (16). The horizontal axis is $x \in [0, L]$.

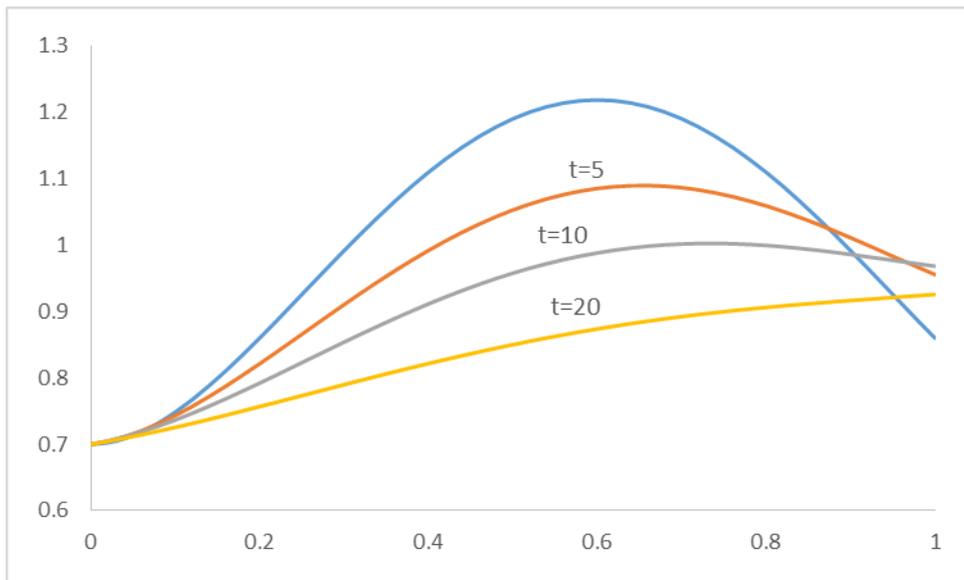

**Fig. 2:** Density ($\rho(t,x)$) profiles for the closed-loop system (1), (14) with (20). The horizontal axis is $x \in [0, L]$.

Figures 5 and 6 show the evolution of the sup norm $\|\rho[t]-\rho^*\|_\infty$ for the initial condition $\rho_0(x) = \rho^* + 4x^2(1.2-x)^2$, $x \in [0,L]$. Figures 5 and 6 exhibit exponential convergence of the sup norm of the state, as predicted by estimate (17). Again, as remarked above, Figures 5 and 6 show that the convergence rate of the closed-loop system (1) with (16) is faster than that the closed-loop system (1), (14) with (20).

As pointed out earlier, the feedback law (20) cannot be used for a set point equal to the critical density. This is not the case for the feedback law (16). Figure 7 shows the density profiles of the closed-loop system (1) with (16) with $\rho^* = \rho_{cr} = 1$ for the initial condition $\rho_0(x) = \rho^* + 4x^2(1.2-x)^2$, $x \in [0,L]$.

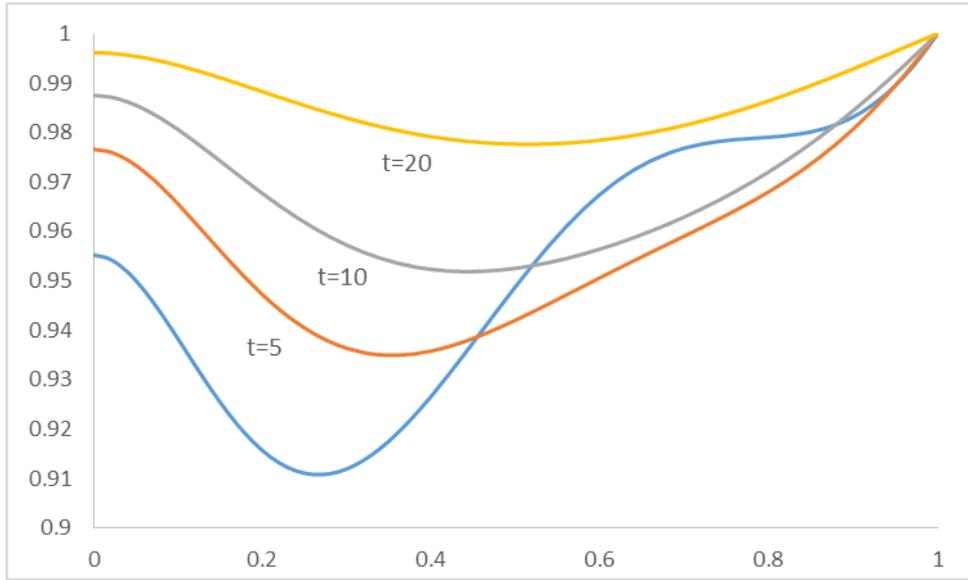

**Fig. 3:** Speed limit ratio ($u(t,x)$) profiles for the closed-loop system (1) with (16). The horizontal axis is $x \in [0,L]$.

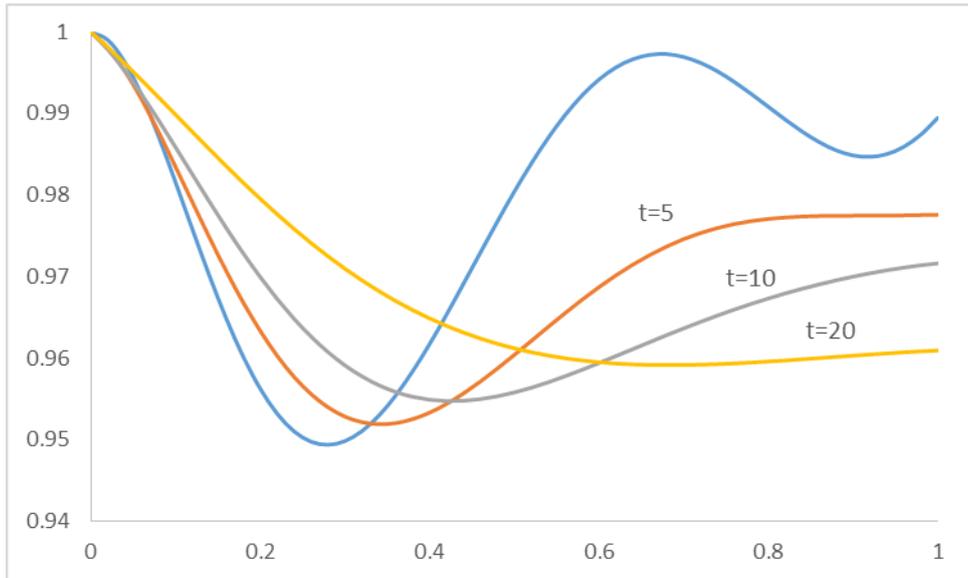

**Fig. 4:** Speed limit ratio ($u(t,x)$) profiles for the closed-loop system (1), (14) with (20). The horizontal axis is $x \in [0,L]$.

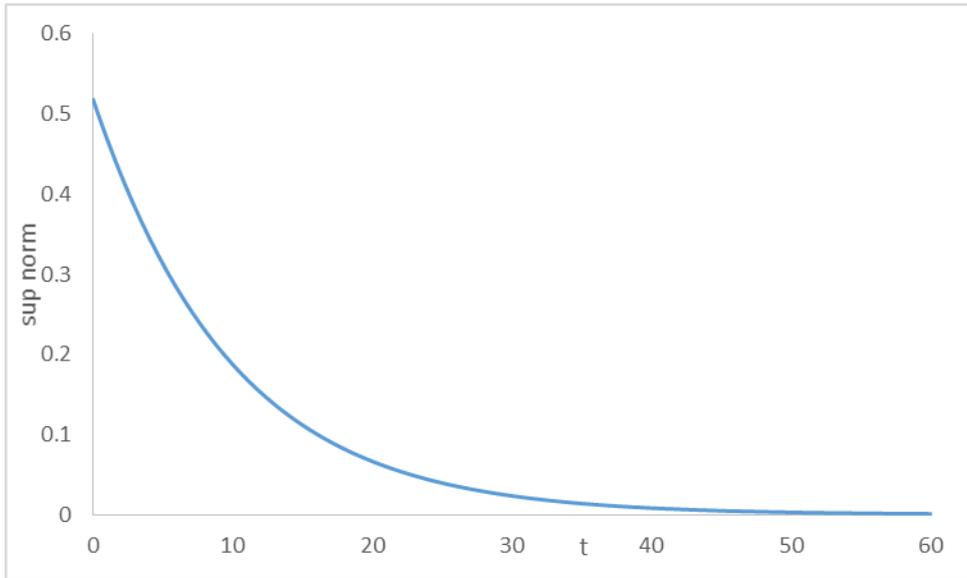

**Fig. 5:** Evolution of the sup norm $\|\rho[t]-\rho^*\|_\infty$ of the solution of the closed-loop system (1) with (16).

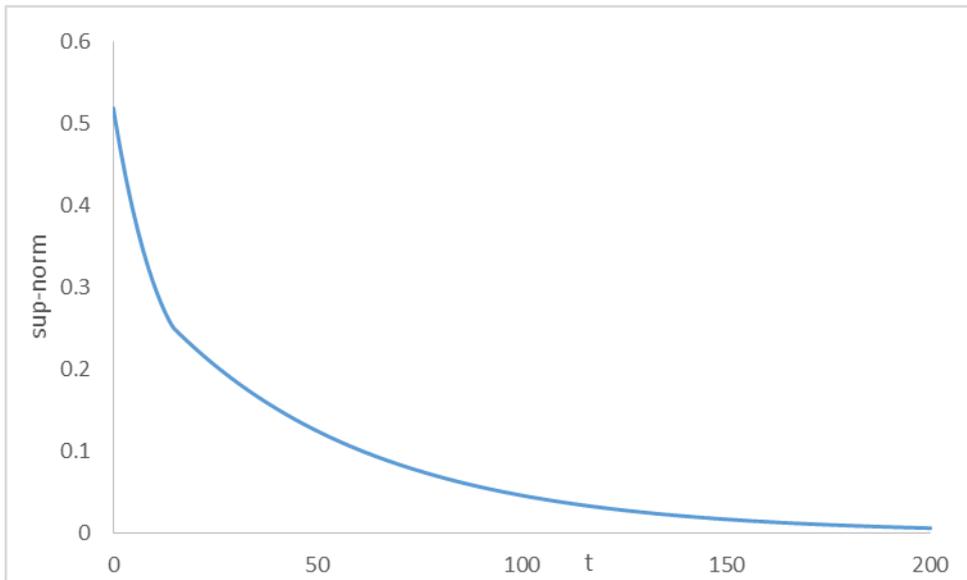

**Fig. 6:** Evolution of the sup norm $\|\rho[t]-\rho^*\|_\infty$ of the solution of the closed-loop system (1), (14) with (20).

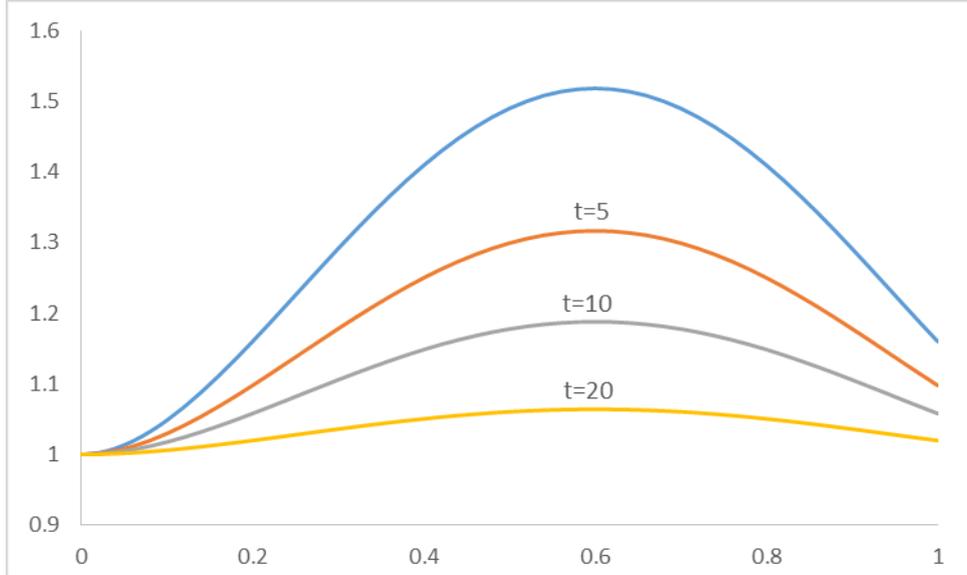

**Fig. 7:** Density ( $\rho(t,x)$ ) profiles for the closed-loop system (1) with (16), $\rho^* = \rho_{cr} = 1$.
The horizontal axis is $x \in [0, L]$.

## 6. Concluding Remarks

The paper provides two different results for the stabilization of a spatially uniform equilibrium profile of a specific scalar conservation law that arises in the study of traffic dynamics under VSL. We have studied two different control problems: the problem with free speed limits at the inlet and the problem with no speed limits at the inlet. For the 1$^{st}$ problem, global asymptotic stabilization was achieved; while for the 2$^{nd}$ problem, regional exponential stabilization was achieved. The solutions of the closed-loop systems were classical solutions, i.e., there are no shocks. The obtained results were illustrated by means of a numerical example.

Future work will address the inhomogeneous case, i.e. consideration of a freeway stretch with spatially inhomogeneous flow-density relationships. Future research may also involve the study of the effect of discretization. In practice, the distributed control input $u(t,x)$ must be kept constant for certain time and space intervals. In other words, the implementation of the feedback controllers (16) or (20) involves discretization both with respect to space and time. The effect of discretization may be important and should be investigated. Another research direction is the incorporation of VSL to 2$^{nd}$ order models. In this way, the velocity dynamics, which are ignored by the 1$^{st}$ order model (1), can be taken into account.

**Acknowledgments:** The authors were supported by ERC funding under the European Union's 7$_{th}$ Framework Programme (FP/2007- 2013)/ ERC Grant Agreement n. [321132], project TRAMAN21.

**References**

[1] Bastin, G., and J.-M. Coron, *Stability and Boundary Stabilization of 1-D Hyperbolic Systems*, Progress in Nonlinear Differential Equations and Their Applications, Springer, 2016.

[2] Blandin, S., X. Litrico, M. L. Delle Monache, B. Piccoli and A. Bayen, "Regularity and Lyapunov Stabilization of Weak Entropy Solutions to Scalar Conservation Laws", *IEEE Transactions on Automatic Control*, 62(5), 2017, 1620 – 1635.

[3] Carlson, R. C., I. Papamichail, M. Papageorgiou and A. Messmer, "Optimal Motorway Traffic Flow Control Involving Variable Speed Limits and Ramp Metering", *Transportation Science*, 44(3), 2010, 238–253.

[4] Carlson, R. C., I. Papamichail, M. Papageorgiou and A. Messmer, "Optimal Mainstream Traffic Flow Control of Large-Scale Motorway Networks", *Transportation Research Part C*, 18, 2010, 193-212.

[5] Carlson, R. C., I. Papamichail and M. Papageorgiou, "Local Feedback-Based Mainstream Traffic Flow Control on Motorways Using Variable Speed Limits", *IEEE Transactions on Intelligent Transportation Systems*, 12, 2011, 1261-1276.

[6] Carlson, R. C., *Mainstream Traffic Flow Control on Motorways*, Ph.D. Thesis, Technical University of Crete, Chania, Greece, 2011.

[7] Coron, J.-M., R. Vazquez, M. Krstic, and G. Bastin. "Local Exponential $H^2$ Stabilization of a 2×2 Quasilinear Hyperbolic System Using Backstepping", *SIAM Journal on Control and Optimization*, 51(4), 2013, 2005–2035.

[8] Coron, J.-M. and G. Bastin, "Dissipative Boundary Conditions for One-Dimensional Quasi-Linear Hyperbolic Systems: Lyapunov Stability for the $C^1$-Norm", *SIAM Journal on Control and Optimization*, 53(4), 2015, 1464-1483.

[9] Cremer, M., *Der Verkehrsfluß auf Schnellstraßen: Modelle, Überwachung*, Regelung, Springer, Berlin, 1979.

[10] Dafermos, C. M., *Hyperbolic Conservation Laws in Continuum Physics*, Dynamical Systems and Differential Equations, Springer, 2010.

[11] Di Meglio, F., R. Vazquez, and M. Krstic, "Stabilization of a System of n+1 Coupled First-Order Hyperbolic Linear PDEs with a Single Boundary Input", *IEEE Transactions on Automatic Control*, 58(12), 2013, 3097-3111.

[12] Goatin, P., S. Göttlich, O. Kolb, "Speed Limit and Ramp Meter Control for Traffic Flow Networks", *Engineering Optimization*, 48 (8), 2016, 1121-1144.

[13] Hegyi, A., *Model Predictive Control for Integrating Traffic Control Measures*, PhD Thesis, TRAIL Thesis Series T2004/2, Delft University of Technology, The Netherlands, 2004.

[14] Hegyi, A., B. De Schutter and J. Hellendoorn, "Optimal Coordination of Variable Speed Limits to Suppress Shock Waves", *IEEE Transactions on Intelligent Transportation Systems*, 6, 2005, 102-112.

[15] Hegyi, A., B. De Schutter, H. Hellendoorn, "Model Predictive Control for Optimal Coordination of Ramp Metering and Variable Speed Limits", *Transportation Research Part C*, 13(4), 2005, 185–209.

[16] Heydecker, B. G. and J. D. Addison, "Analysis and Modelling of Traffic Flow under Variable Speed Limits", *Transportation Research Part C: Emerging Technologies*, 19(3), 2011, 206-217.

[17] Hu, L., F. Di Meglio, R. Vazquez and M. Krstic, "Control of Homodirectional and General Heterodirectional Linear Coupled Hyperbolic PDEs", *IEEE Transactions on Automatic Control*, 61, 2016, 3301–3314.

[18] Karafyllis, I., N. Bekiaris-Liberis and M. Papageorgiou, "Feedback Control of Nonlinear Hyperbolic PDE Systems Inspired by Traffic Flow Models", submitted to *IEEE Transactions on Automatic Control* (see also arXiv:1707.02209 [math.OC]).

[19] Li, T.-T., and W. Libin, *Global Propagation of Regular Nonlinear Hyperbolic Waves*, Birkhauser Boston, 2009.

[20] Li, T. T., *Controllability and Observability for Quasilinear Hyperbolic Systems*, volume 3. Higher Education Press, Beijing, 2009.

[21] Lighthill, M. H., and G. B. Whitham, "On Kinematic Waves II: A Theory of Traffic Flow on Long Crowded Roads", *Proceedings of the Royal Society A*, 229, 1955, 317–345.

[22] Lu, X. Y., T. Z. Qiu, P. Varaiya, R. Horowitz and S. E. Shladover, "Combining Variable Speed Limits with Ramp Metering for Freeway Traffic Control", *Proceedings of the 2010 American Control Conference*, 2010, 2266-2271.

[23] Papageorgiou, M., E. Kosmatopoulos, and I. Papamichail. "Effects of Variable Speed Limits on Motorway Traffic Flow", *Transportation Research Record: Journal of the Transportation Research Board*, 2047, 2008, 37-48.

[24] Prieur, C., and F. Mazenc, "ISS-Lyapunov Functions for Time-Varying Hyperbolic Systems of Balance Laws", *Mathematics of Control, Signals, and Systems*, 24(2), 2012, 111–134.

[25] Prieur, C., J. Winkin, and G. Bastin, "Robust Boundary Control of Systems of Conservation Laws", *Mathematics of Control, Signals, and Systems*, 20(3), 2008, 173–197.

[26] Richards, P. I., "Shock Waves on the Highway", *Operations Research*, 4, 1956, 42–51.

[27] Strub, I., and A. Bayen, "Weak Formulation of Boundary Conditions for Scalar Conservation Laws: An Application to Highway Traffic Modelling", *International Journal of Robust Nonlinear Control*, 16(16), 2006, 733–748.

[28] Yu, H. and M. Krstic, "Traffic Congestion Control on Aw-Rascle-Zhang Model: Full-State Feedback", to appear in the *Proceedings of ACC 2018*.

[29] Yu, H. and M. Krstic, "Adaptive Output Feedback for Aw-Rascle-Zhang Traffic Model in Congested Regime", to appear in the *Proceedings of ACC 2018*.

[30] Zhang, L., and C. Prieur, "Necessary and Sufficient Conditions on the Exponential Stability of Positive Hyperbolic Systems", *IEEE Transactions on Automatic Control*, 62(8), 2017, 3610-3617.